# Four-point bending piezoelectric energy harvester with uniform surface strain toward better energy conversion performance and material usage

Majid Khazaee[1*], John E. Huber[2], Lasse Rosendahl[1], Alireza Rezania[1]

[1] *Department of Energy Technology, Aalborg University, Pontopidanstraede 111, 9220 Aalborg, Denmark.*
[2] *Department of Engineering Science, University of Oxford, Parks Rd, Oxford OX1 3PJ, United Kingdom.*

[*] E-mail: mad@energy.aau.dk (Majid Khazaee)

**Abstract**

Improving the energy conversion efficiency of piezoelectric energy harvesters is of great importance, and one approach is to make more uniform use of the working material by ensuring a uniform strain state. To achieve better performance, this paper presents a four-point bending piezoelectric energy harvester with extensive investigation and modeling to identify the influential parameters. An electromechanical analytical model is presented and verified by experimental data. The frequency-domain method extracts the solutions for a general time-variable force and impact. Four-point bending is compared with the standard cantilever harvesters regarding voltage generation, mechanical strain, and figure of merit. Strain contours are analyzed and interpreted for this innovative approach, and the power generation by the optimal resistance load is studied. Dimensionless parameters are introduced and investigated to find the optimal operating conditions for the four-point bending harvester. Finally, the four-point bending performance and the best figure of merit are discussed with a view to the long-term fatigue life of the harvester. The results show that in the best four-point bending energy conversion conditions; the energy conversion coefficient is more than three times higher than that of typical cantilever energy harvesters. The results also illustrate that the axial strain experienced in a standard cantilever harvester is more than three times higher than that of the four-point bending harvester, suggesting the latter device may have favorable fatigue performance. Overall, the presented piezoelectric harvester has improved energy conversion efficiency and experiences a reduced and uniform surface strain, making it appropriate for high-efficiency energy harvesting systems.

Keywords: Piezoelectric, Energy Harvesting, Harvester configuration, Modeling, Four-Point Bending

# 1. Introduction

Piezoelectric materials, benefiting from electrical-mechanical conversion ability, act as actuators for electrical to mechanical conversion or sensors for mechanical to electrical conversion. The latter can also be used as energy harvesting devices by exploiting the generated electrical energy. Piezoelectric Vibration Energy Harvesting (PVEH) has received much attention from industrial [1] to biomedical [2] applications. By moving toward self-powered electronics by PVEH, electronic devices can be installed in remote areas using the wasted available energy, reducing pollution, production costs, and battery changing costs [3].

In the energy harvesting area, the piezoelectric energy harvester's (PEH) power density, which is the power per unit PEH volume, is significant. The electrical energy generation over the PEH volume should be maximized for high power density. There is a direct link between the mechanical strain and the electrical charge flow in the piezoelectric material. In other words, piezoelectric strain corresponds to electrical energy generation. Efforts have been carried out to enhance power generation by changing strain, such as shape optimization [4] or nonlinear piezoelectric coupling with magnetic forces [5]. On the other hand, in the most used PEH configuration, the clamped-free beam, the strain is non-uniform and small on average over the whole volume [6]. The maximum strain is in the clamped-end region, gradually reducing to zero at the free-end. Thus, a significant harvester volume is ineffective due to the negligible mechanical strain. This low power density issue is specifically problematic where limited volume is available for the PVEH system installation [7].

There have been innovative approaches for improving strain distribution. Modifying the PEH's stiffness by changing the Piezoelectric fiber orientation in the piezoelectric composite harvesters has been investigated [6]. Although this method improved surface strain contour, manufacturing piezoelectric composites with non-zero fiber orientation is challenging, and commercially such samples are unavailable. Changing the location of the piezoelectric material in the cantilevered PEH will change the beam stiffness and affect the surface strain [8]; however, the non-uniform strain in the cantilevered PEH is still present. Strain-engineered material has also been proven effective [9]; however, these laboratory-manufactured materials are not available at large-scale. This work was also tested on a typical cantilevered PEH configuration suffering from non-uniform strain [9]. Nonlinear PEH with the nonlinear substrate structure, auxetic multiple-rotating-cube substrate [10], and perforated substrate [11] under the base excitation has improved PEH performance. Auxetic properties of the piezoelectric material [12] and piezoelectric thickness and poling direction [13] have also been investigated under the cantilevered boundary condition toward improving the PEH performance. Moreover, the PEH shape profile [14] can improve energy generation. The literature shows that most studies have tested or applied innovative PEH design in the cantilevered boundary condition, even in the new bistable [15] and magnetic-piezoelectric coupled [16] systems. Clamped-clamped beams with a center mass have also been investigated in both flat-shape [10] and M-shaped configurations [17]. Therefore, it can be concluded that great attention has been given to the piezoelectric material development or PEH structure modifications for surface strain improvement. In contrast, limited investigation has been given to finding a new setup or boundary condition.

Providing a constant strain over the PEH surface can lead to better material usage and less strain concentration, leading to better energy conversion PEH efficiency. The present study introduces a new way to investigate the uniform surface strain by developing an innovative boundary condition configuration, which is different from the state-of-the-art, focusing on the material/layup improvements. Four-point bending (FPB) is a boundary condition that provides uniform surface strain between the inner clamps; therefore, it is ideal for providing uniform conditions along the length of the PEH. The FPB-PEH has not been investigated, despite its great potential for energy harvesting technologies. The present study indicates improvement in the energy conversion efficiency of piezoelectric energy harvesters with the new FPB configuration and enhances the power density generation by piezoelectric materials. This innovative method with uniform surface strain is also significant because of reduced mechanical strain concentration, better long-term performance, and better usage of the piezoelectric material.

The present study explores an innovative boundary condition for PEHs toward better performance by uniform power generation over the whole PEH surface. Analytical modeling of the PEHs under the FPB condition is presented for a piezoelectric composite material. The rest of the study is categorized as follows. Section 2 first introduces the concept of FPB and its importance. Then, the electromechanically coupled analytical modeling of the FPB harvester is presented with closed-form solutions. Dimensionless parameters and a figure of merit ($FoM$) for the performance evaluation are presented. The results and discussions are presented in section 3. The strain contours, the sensitivity analysis concerning dimensionless parameters, comparisons with the cantilevered-PEH, and the conversion performance of the FPB-PEH are discussed in section 3. Concluding remarks and future directions are presented in section 4. Comparisons of FPB and cantilevered PEHs show that the energy conversion efficiency of the FPB configuration is higher. The experimental tests were carried out to verify the FPB-PEH analytical model. The results also suggest an improved design of the FPB-PEH for better long-term working.

## 2. Methodology

### 2.1. A new energy harvesting configuration: Four-point bending

The cantilevered boundary condition is the most common configuration of PEH [18]; it has non-uniform stress and high-stress regions near the clamp line. These two conditions cause low power generation for low-stress areas, and structural failure becomes imminent for the high-stress parts, making the cantilevered configuration sub-optimal for energy harvesters. Here, numerical results from COMSOL Multiphysics software for the cantilevered beam are compared with the proposed new energy harvesting configuration. Both cantilever and FPB harvesters are analyzed under a static force considering electromechanical coupling effects. The top and bottom of the piezoelectric are covered by electrodes, providing a uniform surface voltage.

Figure 1 (a) shows the von Mises stress of a cantilever PEH with 1 N tip static force, indicating a non-uniform stress distribution over the length; approximately one-third of the beam length is under high stress, and the rest is under low stress. This non-uniform stress distribution leads to non-uniform voltage generation in PEHs; a significant part of the beam does not contribute to the voltage generation, which will be considered a nonoptimal material usage [6].

To overcome the non-uniform stress distribution, the four-point bending (FPB) boundary conditions are proposed. As shown in Figure 1 (b), the PEH is under pin-pin boundary conditions, and the force is applied at two points from the top of the beam. The free diagram of the beam shows a uniform bending moment between the two top forces. Unlike the cantilevered beam, this uniform bending moment leads to a consistent voltage generation in PEHs. Thus, the proposed FPB boundary condition presents a better usage of piezoelectric material for energy generation.

Apart from the uniformity, the stress amplitude near the cantilevered clamp is considerably higher than the FPB harvester. as can be interpreted from Figure 1 (a) and (b). Note that the local stress at single-point forces in FPB depends on the boundary condition; however, considerably lower than the stress near the cantilever harvester. Therefore, by the FPB design, minor stress is applied to the material under the same load, which may benefit the long-term fatigue life of PEHs.

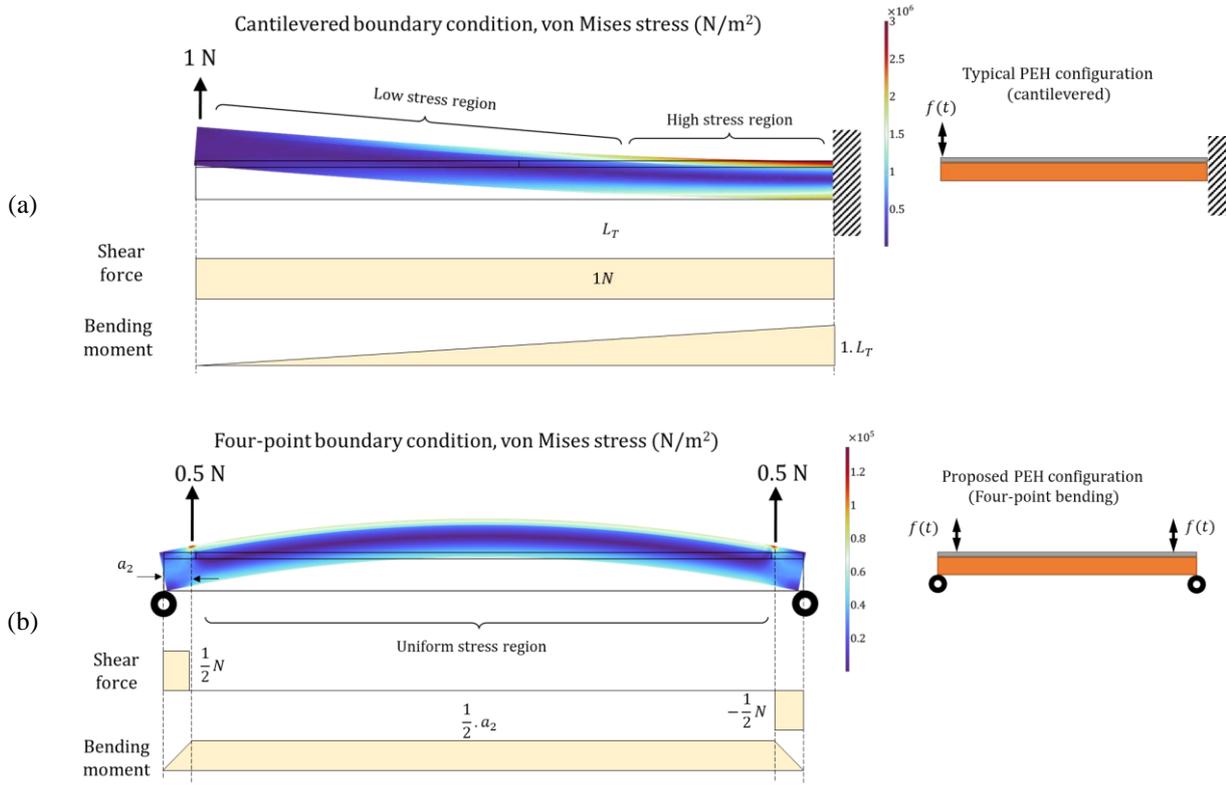

Figure 1. (a) Typical Piezoelectric energy harvester (PEH) and its von Mises stress under 1 N static force, and (b) Four-point bending PEH and its von Mises stress under a total static force of 1 N.

One can design various layouts for the FPB PEH. This study focuses on an impact force on the FPB configuration; thus, considering a simple but practical FPB configuration, a two-clamp design is a valuable way to apply the four-point mechanical boundary conditions in PEHs. The bottom clamp applies the pin-pin free-rotation boundary condition, while the top clamp is for using the load. This two-clamp configuration is shown in Figure 2. The piezoelectric layer partially covers the substrate shim, so the applied force and the boundary clamps act on the substrate shim. These features avoid direct contact force applied to the piezoelectric material, which will be especially beneficial for fragile piezoceramics.

By applying a dynamic force $F(t)$ to the top clamp (TC), this force is divided into two forces $F_B(t)$ at point-B and $F_C(t)$ at point-C, where these forces are calculated from the dynamic equations. For the top clamp, it can be written that

$$\sum F_z = m_{\text{TC}} \frac{\partial^2 z(x,t)}{\partial t^2}\bigg|_{\text{TC}} \rightarrow F(t) - F_B(t) - F_C(t) = m_{\text{TC}} \frac{\partial^2 z(x,t)}{\partial t^2}\bigg|_{\text{TC}} \quad (1)$$

$$\sum M_{\text{CoR}} = I_{\text{CoR}} \alpha \rightarrow F(t) - F_B(t) - F_C(t) = m_{\text{TC}} \frac{\partial^2 z(x,t)}{\partial t^2}\bigg|_{\text{TC}} \qquad (2)$$

where TC stands for the "top clamp", $M$ is the bending moment, CoR is the center of the rotation, $I$ is the rotation inertia, and $\alpha$ is the angular acceleration. According to Eq. (1), a smaller top-clamp mass will be beneficial since the force magnitude applied to the PEH will then be greater. More importantly, the top-clamp mass effect should be considered for modeling purposes. In Eq. (2), $\alpha \neq 0$ means that the top clamp rotates during the force application, which is undesirable. For having $M_{CoR} = 0$, the rotation inertia should be zero, or the clamp conditions should be symmetric. Therefore, the symmetric four-point bending energy harvesters were used in our experiments to minimize rotation, while the general unsymmetric case is considered for the mathematical modeling.

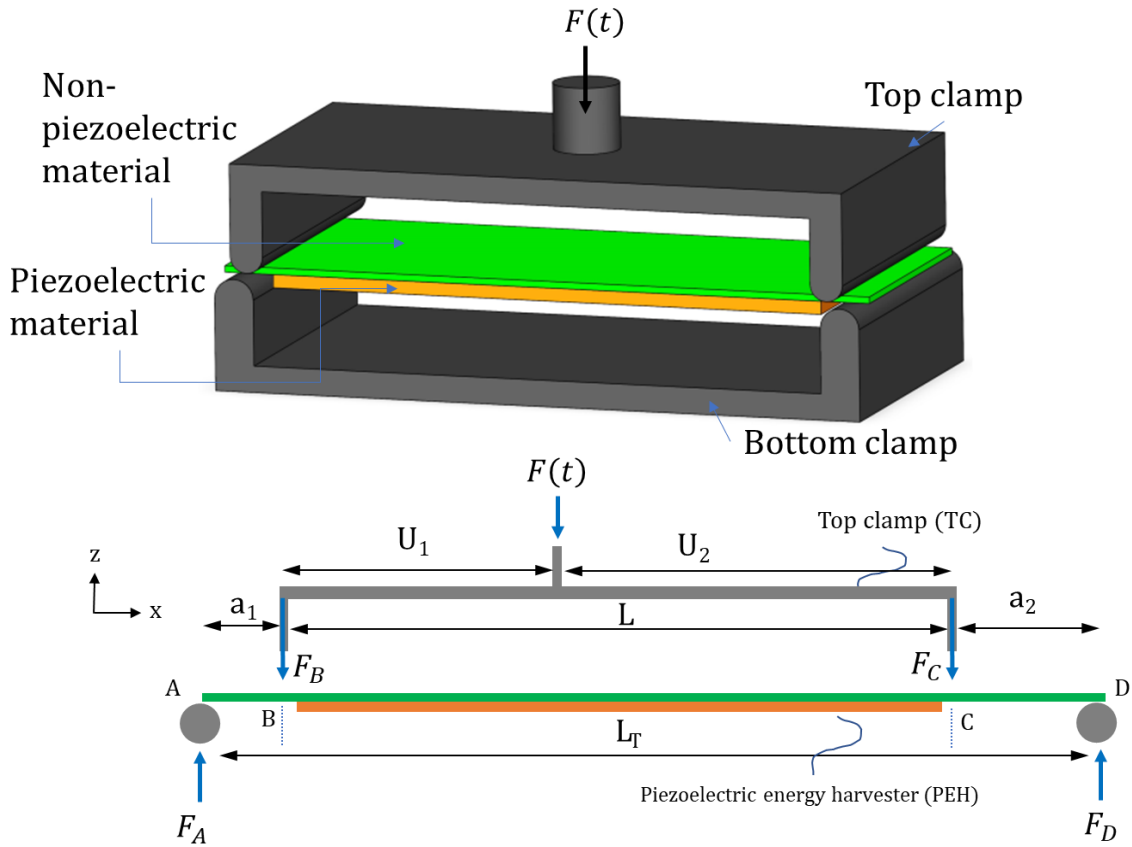

Figure 2. A two-clamp configuration for the FPB harvester.

### 2.2. Modeling the piezoelectric beam under the four-point bending condition

### 2.2.1. Mechanical vibration equation

This section presents the electromechanical-coupled vibration equation in the four-point bending (FPB) boundary condition. As a general case, the piezoelectric element is a multi-layered beam with active embedded piezoelectric layers. Thus, the vibration equation is first presented for a broad multi-layered beam, and later for cases of unimorph and bimorph, the detailed vibration equation is then examined.

The FPB condition is shown in Figure 3 (a). The PEH rests on two pins. The beam is subjected to perpendicular time-dependent external forces, $F_B$ and $F_C$. The piezoelectric beam in this boundary condition is a 31-mode energy harvester, where the poling is in the thickness direction (3 direction), and the mechanical strain of interest is in the length direction (1 direction). The two output PEH wires are connected to a purely resistive electrical load, $R_L$.

The PEH's cross-section for multi-layered, unimorph, and bimorph layup is shown in Figure 3 (b). Note the neutral axis location, which affects the beam's structural properties. The neutral axis is located at the substrate middle line for symmetric layups. The neutral axis location for nonsymmetric layups shall be calculated numerically from the material and thickness variables.

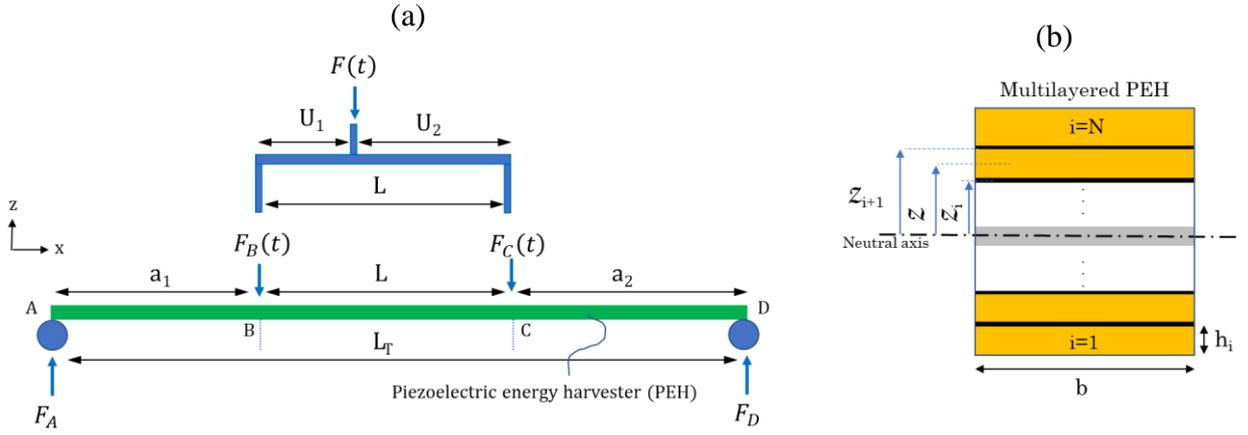

Figure 3. (a) Four-point bending piezoelectric energy harvester (FPB-PEH) and (b) side-view of the harvester with multi-layered, unimorph, and bimorph configurations.

The total force distribution along the beam can be expressed by:

$$Q(x,t) = F_A \cdot \delta(0) - F_B \cdot \delta(a_1) - F_C \cdot \delta(L_T - a_2) + F_D \cdot \delta(L_T) \tag{3}$$

The force distribution can be further simplified using the static relations between the $F_A$, $F_B$, $F_C$, and $F_D$ components, as given by:

$$Q(x,t) = \left(\left(1 - \frac{a_1 + U_1}{L_T}\right)\delta(0) - \left(1 - \frac{U_1}{L}\right)\delta(a_1) - \left(\frac{U_1}{L}\right)\delta(L_T - a_2) \right. \\ \left. + \left(\frac{a_1 + U_1}{L_T}\right)\delta(L_T)\right) F(t) \tag{4}$$

The general equation of motions for a bending beam based on the Euler-Bernoulli beam theory can be expressed as [19]:

$$\frac{\partial^2 M(x,t)}{\partial x^2} + C_a \frac{\partial w(x,t)}{\partial t} + m^* \frac{\partial^2 w(x,t)}{\partial t^2} = Q(x,t) \tag{5}$$

where $m^*$ is the harvester mass per length expressed by $m^* = \sum_{i=1}^{N} bh_i\rho_i$.

Internal bending moment, $M(x,t)$, for a multi-layered harvester is calculated by the axial stress, as given by

$$M(x,t) = -b\left(\sum_{i=1}^{N} \int_{Z_i}^{Z_{i+1}} T_{xx} Z dz\right) \tag{6}$$

wherein $Z$ is the distance from the neutral axis and $T_{xx}$ is the axial stress.

The axial stress ($T_{xx}$) is related to the beam deformation $w(x,t)$ by the constitutive and strain-deformation equations. The linear constitutive equations for the piezoelectric and substrate elastic layers are respectively given by [20]

$$T_{xx}^p = \bar{c}_{11}^E \varepsilon_{xx} - \bar{e}_{31} E_z \tag{7.a}$$

$$T_{xx}^s = Y_s \varepsilon_{xx} \tag{7.b}$$

where $\varepsilon_{xx}$ is the axial strain, $E_z$ is the electric field, and $\bar{e}_{31}$ is the piezoelectric coupling coefficient. The terms in Eq. (7.a) and Eq. (7.b) shall be further elaborated as follows:

- $\boldsymbol{\varepsilon_{xx}}$: In the linear framework of a pure bending beam, the axial strain at a certain level $z$ from the neutral axis can be obtained from the beam curvature, as given by

$$\varepsilon_{xx}(x,z,t) = -z\frac{\partial^2 w(x,t)}{\partial x^2} \tag{8}$$

- $\boldsymbol{E_z}$: Under a constant electrical field, $E_z$ can be related to the voltage and piezoelectric layer's thickness. The electrical field for a single piezoelectric layer is $E_z = -\frac{V_R}{h_p}$. In the multilayer configuration, note the reverse electric field signs for the piezoelectric layers above and below the neutral axes. For series connection $E_z = -\frac{V_R}{2h_p}$ and for parallel connection $E_z = -\frac{V_R}{h_p}$ [21].

Substituting Eq. (7.a), Eq. (7.b), and Eq. (8) into Eq. (6) gives the general form of the beam bending stiffness as:

$$\mathrm{M}(x,t) = YI\frac{\partial^2 \mathrm{w}(x,t)}{\partial x^2} + \mathcal{P}V_R(t)\big(\mathcal{H}(x-x_i) - \mathcal{H}(x-x_f)\big) \qquad (9)$$

where $YI$ is the bending stiffness, $\mathcal{P}$ is the electromechanical coupling factor, and $\big(\mathcal{H}(x-x_i) - \mathcal{H}(x-x_f)\big)$ ensures that the piezoelectric layer spans from $x_i$ to the $x_f$. In our case studies, $x_i = 0$ and $x_f = L_T$.

Substituting Eq. (9) into the beam vibration equation, Eq. (5), the differential equations of motion can be expressed as a function of the relative beam deflection, as given by: (note that $\frac{d\mathcal{H}(x)}{dx} = \delta(x)$)

$$YI\frac{\partial^4 \mathrm{w}(x,t)}{\partial x^4} + C_a\frac{\partial \mathrm{w}(x,t)}{\partial t} + m^*\frac{\partial^2 \mathrm{w}(x,t)}{\partial x^2} + \mathcal{P}V_R(t)\left(\frac{d\delta(x-x_i)}{dx} - \frac{d\delta(x-x_f)}{dx}\right) = Q(x,t) \qquad (10)$$

The four-point bending (Figure 3) boundary conditions for Eq. (10) are expressed with,

At $x = 0 \rightarrow \mathrm{w}(x=0,t) = 0$ and $\mathrm{M}(x=0,t)=0$ \hfill (11.a)

At $x = L_T \rightarrow \mathrm{w}(x=L_T,t) = 0$ and $\mathrm{M}(x=L_T,t) = 0$ \hfill (11.b)

From the bending moment definition in Eq. (9), the bending moment boundary condition can also be given by:

$$\left.\frac{\partial^2 \mathrm{w}(x,t)}{\partial x^2}\right|_{x=0,L_T} = -\frac{\mathcal{P}V_R(t)}{YI} \qquad (11.c)$$

$N.m^2 = \mathcal{P}.V$

In practice, $|\mathcal{P}|$ is in the range of $10^{-5}$ to $2\times 10^{-5}$ $C$, and $YI$ is in the range of $5\times 10^{-3}$ to $2\times 10^{-2}$ $N.m^2$ for the unimorph and bimorph, respectively. These values suggest that $\left|\frac{\mathcal{P}V_R(t)}{YI}\right|$ is negligible, at about 1% for a 10 V voltage generation ($V_R=10$ V) for the unimorph and bimorph, respectively. Therefore, approximately, the bending moment boundary conditions can be simplified into

$$\left.\frac{\partial^2 \mathrm{w}(x,t)}{\partial x^2}\right|_{x=0 \text{ and } L_T} = 0 \qquad (11.d)$$

The beam deflection is expanded as an infinite series of mode shapes for solving the differential equation, Eq. (10), as given by:

$$w(x,t) = \sum_{n=1}^{\infty} \Pi_n(t)\, \phi_n(x) \tag{12}$$

where mode shapes, $\phi_n(x)$, should satisfy the boundary conditions, and $\Pi_n(t)$ is the time-dependent mode contribution. For the boundary conditions in Eq. (11.a) and (11.b), the mode shapes are given by [22]:

$$\phi_n(x) = C_n \sin\left(n\pi \frac{x}{L_T}\right) \tag{13}$$

Note that the modal orthogonality is satisfied by the mode shapes, i.e., $\int_0^{L_T} \phi_m(x)\phi_n(x)dx = \delta_{mn}$. The constant $C_n$ is a normalization coefficient for satisfying modal orthogonality. For the symmetric boundary condition, the deflection at B and C are the same ($w]_{@B} = w]_{@C}$) and therefore, only the odd selections of $n$ are valid.

Substituting Eq. (12) into the differential Eq. (10) gives:

$$m^* \sum_{n=1}^{\infty} \phi_n(x)\ddot{\Pi}_n(t) + C_a \sum_{n=1}^{\infty} \phi_n(x)\dot{\Pi}_n(t) + YI\left(\frac{n\pi}{L_T}\right)^4 \sum_{n=1}^{\infty} \phi_n(x)\Pi_n(t) \\ + \mathcal{P}V_R(t)\left(\frac{d\delta(x-x_i)}{dx} - \frac{d\delta(x-x_f)}{dx}\right) = Q(x,t) \tag{14}$$

Multiplying Eq. (14) by $\frac{\phi_m(x)}{m^*}$ from the left side and integrating from 0 to $L_T$, Eq. (14) will simplifies to

$$\ddot{\Pi}_n(t) + 2\zeta_n \omega_n \dot{\Pi}_n(t) + \omega_n^2 \Pi_n(t) + \gamma_n V_R(t) = \frac{1}{m^*}\sigma_n F(t) \tag{15}$$

where $\zeta_n$ is the modal damping coefficient, $\omega_n$ is the natural frequency and $\gamma_n$ is the modal electromechanical coupling factor, which is given by:

$$\omega_n^2 = \left(\frac{n^2\pi^2}{L_T^2}\right)^2 \frac{YI}{m^*} \tag{16.a}$$

$$\zeta_n = \frac{C_a}{2m^*\omega_n} \tag{16.b}$$

$$\gamma_n = \frac{\mathcal{P}}{m^*}\left(\frac{d\phi_n(x)}{dx}\bigg|_{x=x_i}^{x=x_f}\right) \tag{16.c}$$

For the force coefficient, there will be only two non-zero elements in $Q(x,t)$, because $\phi_n(x=0) = \phi_n(x=L_T) = 0$; therefore, $Q(x,t) = -\left(\left(1-\frac{U_1}{L}\right)\delta(a_1) + \left(\frac{U_1}{L}\right)\delta(L_T - a_2)\right)F(t)$.

Consequently, the force coupling coefficient becomes:

$$\sigma_n = -\left(\left(1-\frac{U_1}{L}\right)\cdot\phi_n|_{x=a_1} + \frac{U_1}{L}\cdot\phi_n|_{x=L_T-a_2}\right) \tag{17}$$

### 2.2.2. The electrical equation for voltage output

The electrical equation for a single Piezoceramic layer is derived. Then, the equivalent electrical equation can be obtained according to the electrical connections of multiple piezoceramic layers (series or parallel).

According to Gauss' law,

$$\frac{d}{dt}\left(\iint_A \vec{D}\cdot\vec{n}\,dA\right) = \frac{d}{dt}\left(\iint_A D_3\,dA\right) = \frac{V_R(t)}{R_L} \tag{18}$$

where $\vec{D}$ is the electrical displacement vector, $A$ is the surface of one electrode, and $\vec{n}$ is the unit normal vector to the $A$.

The piezoelectric constitutive equation for the electrical displacement states that [20]

$$D_3 = \bar{e}_{31}\varepsilon_{xx} + \bar{\epsilon}_{33}E_z \tag{19}$$

wherein $\bar{\epsilon}_{33}$, the permittivity constant at constant strain is obtained from the permittivity constant ($\epsilon_{33}^T$) at constant stress by $\epsilon_{33}^T = \bar{\epsilon}_{33} - e_{31}^2/\bar{c}_{11}^p$.

For a constant electric field across the piezoelectric thickness, the electric field in the $z$-direction is $E_z = -\frac{V_R}{h_p}$.

Substituting $\varepsilon_{xx}$ in Eq. (8) into Eq. (18), Gauss' law yields:

$$\frac{\bar{\epsilon}_{33}b(x_f - x_i)}{h_p}\frac{dV_R(t)}{dt} + \frac{V_R(t)}{R_L} = -Z_p\bar{e}_{31}b\int_{x_i}^{x_f}\frac{\partial^3 w(x,t)}{\partial t\,\partial x^2}dx \tag{20}$$

Note that $Z_p$ is the distance of the piezoelectric mid-plane to the neutral axis. Analogous to the load-capacitor electrical circuit equations, the piezoelectric capacitance coefficient from Eq. (20) is the multiplier of $\frac{dV_R(t)}{dt}$.

The mechanical coupling term of Eq. (20) is simplified using the modal expansion assumption in Eq. (12), as given by

$$\int_{x_i}^{x_f} \frac{\partial^3 w(x,t)}{\partial t\, \partial x^2} dx = \sum_{n=1}^{\infty} \left( \int_{x_i}^{x_f} \frac{d^2 \phi_n(x)}{dx^2} dx \right) \dot{\Pi}_n(t) = \sum_{n=1}^{\infty} \left( \frac{d\phi_n(x)}{dx} \bigg|_{x_i}^{x_f} \right) \dot{\Pi}_n(t) \qquad (21)$$

Therefore, the electrical equation is simplified to:

$$C_P \frac{dV_R(t)}{dt} + \frac{V_R(t)}{R_L} = I_p(t) \qquad (22)$$

wherein the piezoelectric layer capacitance ($C_P$), the modal coupling factor ($\Lambda_n$), and the piezoelectric current source are respectively given by Eq. (23), Eq. (24), and Eq. (25).

$$C_P = \frac{\bar{\epsilon}_{33} b (x_f - x_i)}{h_p} \qquad (23)$$

$$\Lambda_n = -\mathcal{Z}_p \bar{e}_{31} b \left( \frac{d\phi_n(x)}{dx} \bigg|_{x_i}^{x_f} \right) \qquad (24)$$

$$I_p(t) = \sum_{n=1}^{\infty} \Lambda_n \dot{\Pi}_n(t) \qquad (25)$$

Eq. (22) is derived for one piezoelectric layer. An equivalent electrical circuit can be considered for multiple piezoelectric layers.

### 2.2.3. Special Cases of unimorph and bimorph harvesters

The previous subsections presented the analysis of multi-layered piezoelectric energy harvesters. In many practical energy harvesters, a single piezoelectric layer (unimorph) or double piezoelectric layers (bimorph) are employed; thus, these two cases are separated here.

Figure 4 (a) shows the piezoelectric unimorph and bimorph harvesters' side-view with the thickness properties. The neutral axis location is shown for them. Moreover, the electrical circuits for unimorph and bimorph (series and parallel connections) are shown in Figure 4 (b). These different electrical circuits can be transformed into an equivalent circuit, allowing a unified electrical equation.

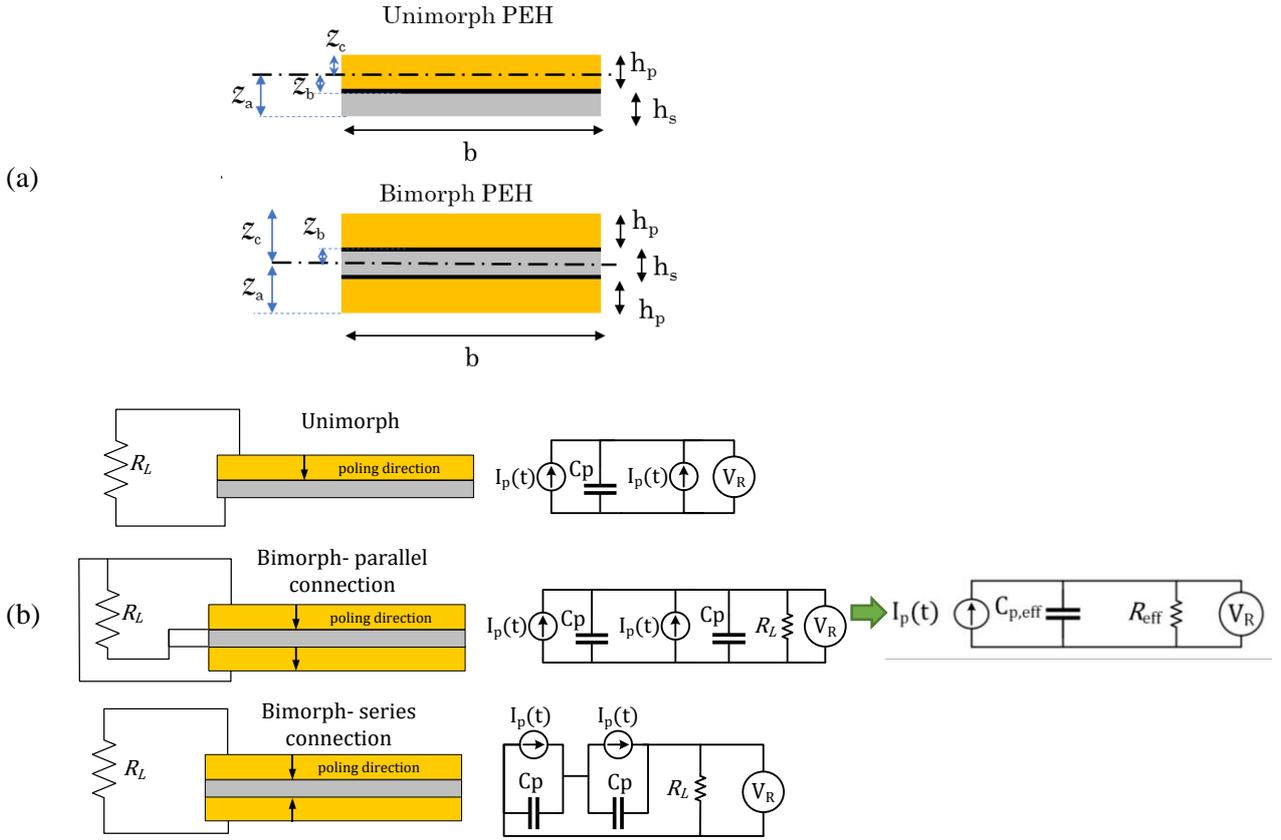

Figure 4. (a) side-view of the unimorph and bimorph PEHs with (b) their equivalent circuits.

The internal bending moment for the unimorph and bimorph is given by

Unimorph: $\quad M(x,t) = -b \left( \int_{-Z_a}^{-Z_b} T_{xx}^s z\, dz + \int_{-Z_b}^{Z_c} T_{xx}^p z\, dz \right)$ (26.a)

Bimorph: $\quad M(x,t) = -b \left( \int_{-\frac{h_s}{2}-h_p}^{-\frac{h_s}{2}} T_{xx}^p z\, dz + \int_{-\frac{h_s}{2}}^{\frac{h_s}{2}} T_{xx}^s z\, dz + \int_{\frac{h_s}{2}}^{\frac{h_s}{2}+h_p} T_{xx}^p z\, dz \right)$ (26.b)

The bending stiffness and coupling factors for the unimorph and bimorph beams are simplified into:

Unimorph:
$$YI = {b}/{3} \left[ Y_s(Z_b^3 - Z_a^3) + \bar{c}_{11}^E(Z_c^3 - Z_b^3) \right].$$
$$\mathcal{P} = -\frac{\bar{e}_{31} b}{2h_p}[Z_c^2 - Z_b^2]$$
(27.a)

Bimorph: $YI = {}^{2b}/_3 \left[ Y_s {}^{h_s^3}/_8 + \bar{c}_{11}^E \left( \left( h_p + {}^{h_s}/_2 \right)^3 - {}^{h_s^3}/_8 \right) \right]$

— Series connection $\quad \mathcal{P} = \dfrac{\bar{e}_{31} b}{2 h_p} \left[ {}^{h_s^2}/_4 - \left( h_p + {}^{h_s}/_2 \right)^2 \right]$ (27.b)

— Parallel connection $\quad \mathcal{P} = \dfrac{\bar{e}_{31} b}{h_p} \left[ {}^{h_s^2}/_4 - \left( h_p + {}^{h_s}/_2 \right)^2 \right]$

Regarding the electrical equation simplification for unimorph and bimorph, the distance of the piezoelectric center line to the neutral axis is $\mathcal{Z}_p = (h_s + h_p)/2$, and thus the coupling term in the electrical equation is:

$$\Lambda_n = -\dfrac{(h_s + h_p)}{2} \bar{e}_{31} b \left( \dfrac{d\phi_n(x)}{dx} \bigg|_{x_i}^{x_f} \right) \quad (28)$$

The equivalent circuit elements for the unimorph and bimorph harvesters, as shown in Figure 4 (b), are given by:

Unimorph $\qquad\qquad\qquad\qquad C_{P,\text{eff}} = C_P,\ R_{\text{eff}} = R_L \qquad$ (29.a)

Bimorph in series connection $\qquad C_{P,\text{eff}} = C_P,\ R_{\text{eff}} = R_L \qquad$ (29.b)

Bimorph in parallel connection $\qquad C_{P,\text{eff}} = C_P,\ R_{\text{eff}} = 2 R_L \qquad$ (29.c)

### 2.2.4. Closed-form solutions for the mechanical and electrical responses

Frequency domain analysis is employed for the closed-form solutions under a general load. In the frequency domain analysis, it is considered that the external force is represented by a series of harmonic functions using the Fourier Transform (FT):

$$\widehat{F}(\omega) = \int_{-\infty}^{\infty} F(t)\, e^{-j\omega t} dt \quad (30)$$

wherein $\widehat{F}(\omega)$ is the external force FT.

In the linear framework, the mechanical displacement $\Pi(t)$ and the output voltage $V_R(t)$ can be extracted by summing up the output from each harmonic $\omega$, starting from $-\infty$ to $+\infty$. Individual harmonic components of mechanical displacement and voltage are denoted by $\bar{\eta}_\omega$ and $\bar{V}_{R,\omega}$, respectively. Therefore, the overall mechanical vibration and electrical equations can be obtained by integration, as

$$\Pi(t) = \int_{-\infty}^{\infty} \bar{\Pi}_\omega e^{j\omega t} d\omega \cong \sum_{r=-\infty}^{\infty} \bar{\Pi}_{\omega_r} \cdot e^{j\omega_r t} \Delta\omega_r \qquad (31.a)$$

$$V_R(t) = \int_{-\infty}^{\infty} \bar{V}_{R,\omega} e^{j\omega t} d\omega \cong \sum_{r=-\infty}^{\infty} \bar{V}_{R,\omega_r} \cdot e^{j\omega_r t} \Delta\omega_r \qquad (31.b)$$

wherein $\bar{\Pi}_{\omega_r}$ is the mechanical response and $\bar{V}_{R,\omega_r}$ is the piezoelectric voltage response due to a nominal harmonic excitation with $\omega_r$ frequency. Note that the over-bar indicates the magnitude.

For obtaining $\bar{\Pi}_{\omega_r}$ and $\bar{V}_{R,\omega_r}$, the harmonic solution of the piezoelectric energy harvester differential equations should be obtained. The governing differential equations are:

$$\ddot{\Pi}_{n,\omega}(t) + 2\zeta_n\omega_n\dot{\Pi}_{n,\omega}(t) + \omega_n^2\Pi_{n,\omega}(t) + \gamma_n V_{R,\omega}(t) = \frac{\sigma_n \hat{F}(\omega)}{2\pi m^*} e^{j\omega t} \qquad (32.a)$$

$$C_{P,\text{eff}} \frac{dV_{R,\omega}(t)}{dt} + \frac{V_{R,\omega}(t)}{R_{\text{eff}}} = \sum_{n=1}^{\infty} \Lambda_n \dot{\Pi}_{n,\omega}(t) \qquad (32.b)$$

By substituting series representations of the mechanical and the output voltage from (31), the steady-state relationships can be given by,

$$(\omega_n^2 - \omega^2 + j2\zeta_n\omega_n\omega)\bar{\Pi}_{n,\omega} + \gamma_n \bar{V}_{R,\omega} = \frac{\sigma_n \hat{F}(\omega)}{2\pi m^*} \qquad (33.a)$$

$$\left(\frac{1}{R_{\text{eff}}} + jC_{P,\text{eff}}\omega\right)\bar{V}_{R,\omega} = \sum_{n=1}^{\infty} j\omega\Lambda_n \bar{\Pi}_{n,\omega} \qquad (33.b)$$

Eliminating the mechanical response $\bar{\eta}_{n,r}$ between Eq. (33.a) and (33.b), the output voltage can be expressed as,

$$\bar{V}_{R,\omega} = \frac{\hat{F}(\omega)}{2\pi} \Psi(\omega, R_{\text{eff}}) \qquad (34)$$

wherein $\Psi(\omega, R_{\text{eff}})$ is the energy conversion term defined by

$$\Psi(\omega, R_{\text{eff}}) = \frac{\frac{1}{m^*} j\omega \sum_{n=1}^{\infty} \Lambda_n \sigma_n \alpha_{n,\omega}}{\frac{1}{R_{\text{eff}}} + j\omega C_{P,\text{eff}} + j\omega \sum_{n=1}^{\infty} \Lambda_n \gamma_n \alpha_{n,\omega}} \qquad (35)$$

This energy conversion term is frequency- and load-dependent, in addition to depending on the material and geometrical properties. Furthermore, the mechanical frequency response function (FRF), $\alpha_{n,\omega}$, is defined by:

$$\alpha_{n,\omega} = \frac{1}{\omega_n^2 - \omega^2 + j(2\zeta_n \omega_n \omega)} \qquad (36)$$

Substituting Eq. (34) into Eq. (33.a), the mechanical response can be obtained, as expressed by

$$\overline{\Pi}_{n,\omega} = \frac{\hat{F}(\omega)}{2\pi} \alpha_n(\omega) \left[ \frac{\sigma_n}{m^*} - \gamma_n \Psi(\omega, R_{\text{eff}}) \right] \qquad (37)$$

Because of the FRF properties at the resonance $\alpha_n(\omega_r = \omega_n)" \alpha_n(\omega_r \neq \pm\omega_n)$ [23], the main contributions in the sampled frequency $\omega_r$ are those at the natural frequencies of the harvester. This so called the modal assumption, implies that the contribution of mechanical and voltage responses of the piezoelectric beam is prominent at the modal frequencies. Considering the harvester's natural frequencies, the $\alpha_n$ and $\Psi$ can be simplified by the modal assumption, as

$$\alpha_n(\omega = \omega_n) = \frac{1}{j(2\zeta_n \omega_n^2)} \qquad (38.a)$$

$$\Psi(\omega = \omega_n, R_{\text{eff}}) = \frac{1}{m^*} \sum_{n=1,\dots,\infty} \frac{\sum_{n=1}^{\infty} \frac{\Lambda_n \sigma_n}{2\zeta_n \omega_n}}{\frac{1}{R_{\text{eff}}} + jC_{P,\text{eff}}\omega_n + \sum_{n=1}^{\infty} \frac{\Lambda_n \gamma_n}{2\zeta_n \omega_n}} \qquad (38.b)$$

In the impact case, the external force takes the form of a Dirac impulse function $F(t) = F_0 \delta(t - t_0)$, and the Fourier Transform is $\hat{F}(\omega) = F_0 e^{-j\omega t_0}$. By applying the modal assumption for the impact force, the force-normalized voltage output and the maximum beam deflection are given by, respectively

$$\frac{V_R(t)}{F_0} = \frac{1}{2\pi} \sum_{n=1}^{N} \Psi(\omega_n, R_{\text{eff}}) \tag{39}$$

$$\frac{w(x,t)|_{\max}}{F_0} \sim \frac{1}{2\pi} \sum_{n=1}^{N} \alpha_n(\omega_n) \left[\frac{\sigma_n}{m^*} - \gamma_n \Psi(\omega_n, R_{\text{eff}})\right] \tag{40}$$

Wherein $N$ is the maximum mode number in consideration.

The force-normalized axial stress from Eq. (8) can be calculated using this beam deflection function, as given by,

$$\frac{\varepsilon_{xx}(x,z,t)|_{\max}}{F_0} = -z \frac{1}{F_0} \frac{\partial^2 w(x,t)}{\partial x^2} = \frac{Z_{\max} \pi}{2L_T^2} \sum_{n=1}^{N} \alpha_n(\omega_n) \left[\frac{\sigma_n}{m^*} - \gamma_n \Psi(\omega_n, R_{\text{eff}})\right] \tag{41}$$

Two parameters are defined for assessing the FPB harvester performance: $\lambda$ (the energy conversion index) and $FoM$ (Figure of Merit). As a metric for demonstrating the voltage generation per unit strain, $\lambda$ is defined as the ratio of open-circuit voltage to the maximum strain, as expressed by

$$\lambda = \frac{V_{OC}}{\max(\varepsilon_{xx})} \tag{42}$$

A higher $\lambda$ means higher voltage generation for a given strain on the PEH. Greater $\lambda$ is advantageous for energy harvesting applications. Nevertheless, $\lambda$ only considers the harvester's performance, not the voltage generated. Thus, the $FoM$ is defined as the product of $\lambda$ and V; in other words, it is the product of voltage generation and energy conversion. For the sake of unit compatibility, the $FoM$ is defined by the square root of the product, as given by

$$FoM = \sqrt{\lambda V_{OC}} \tag{43}$$

## 3. Results and discussion

### 3.1. The reference energy harvester

The result section employs a reference energy harvester in the four-point bending (FPB) boundary condition. Figure 5 shows this reference piezoelectric energy harvester sample and FPB condition with its dimensions. Moreover, the harvesters' dimensionless parameters are defined to provide dimension-independent results. The energy harvester consists of a Macro-Fiber-Composite (MFC), a copper substrate shim, and a double-layer tape as the bonding layer. The impact is applied for two upper clamp spans, with a variable span, e.g., $L$. The reference energy harvester force span is $L=90$ mm ($\kappa=0.83$). For $L>85$ mm, the impact force is not applied to the piezoelectric composite material, so there will be less local impact strain on the MFC and the bending moment over the MFC is constant. The reason for selecting MFC as a harvesting element is that the MFC is a flexible material, not fragile, with reasonably good conversion efficiency. The choice of the bonding layer is due to this material's low material damping, as previously shown by the authors in Ref. [24], resulting in higher power output.

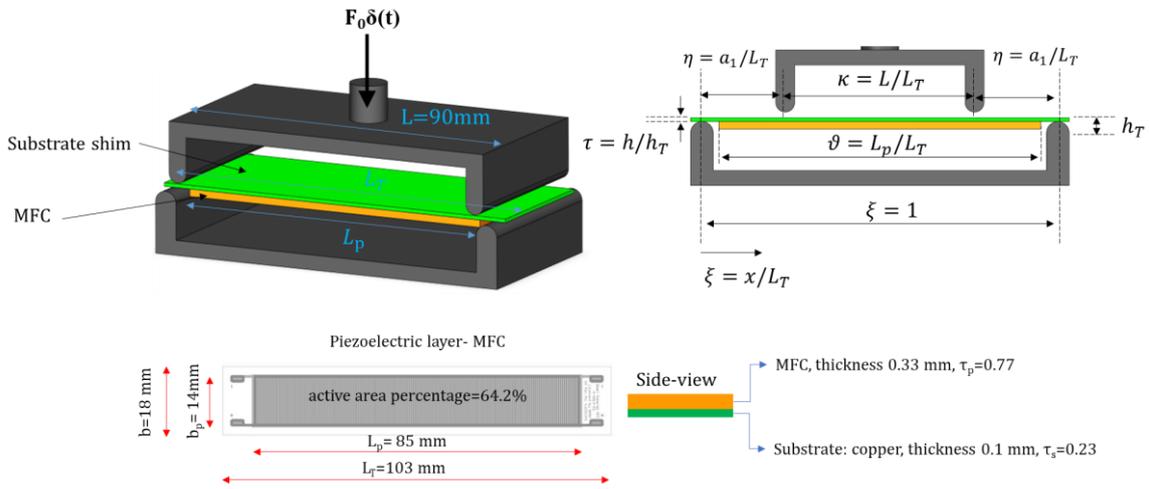

Figure 5. The reference four-point bending (FPB) Piezoelectric energy harvester and dimensionless parameters

## 3.2. The beam model verification by experimental and Finite Element results

The presented FPB model is verified against experimental data in this subsection. The PEH rests on a two-round-end clip; therefore, the beam's vertical motion is restrained, but the rotation is free. A two-round-end clip from the top applies the impact force (see Figure 6). The upper clamp clip contacts the PEH only when providing the impact force, so the upper clip weight does not act on the PEH. A force transducer measures the hammer's impact, and the National Instrument voltage input module measures the piezoelectric output voltage. Figure 6 (a) shows the hammer experimental verification setup. Two FPB cases have been applied to the unimorph reference harvester, shown in Figure 6 (b).

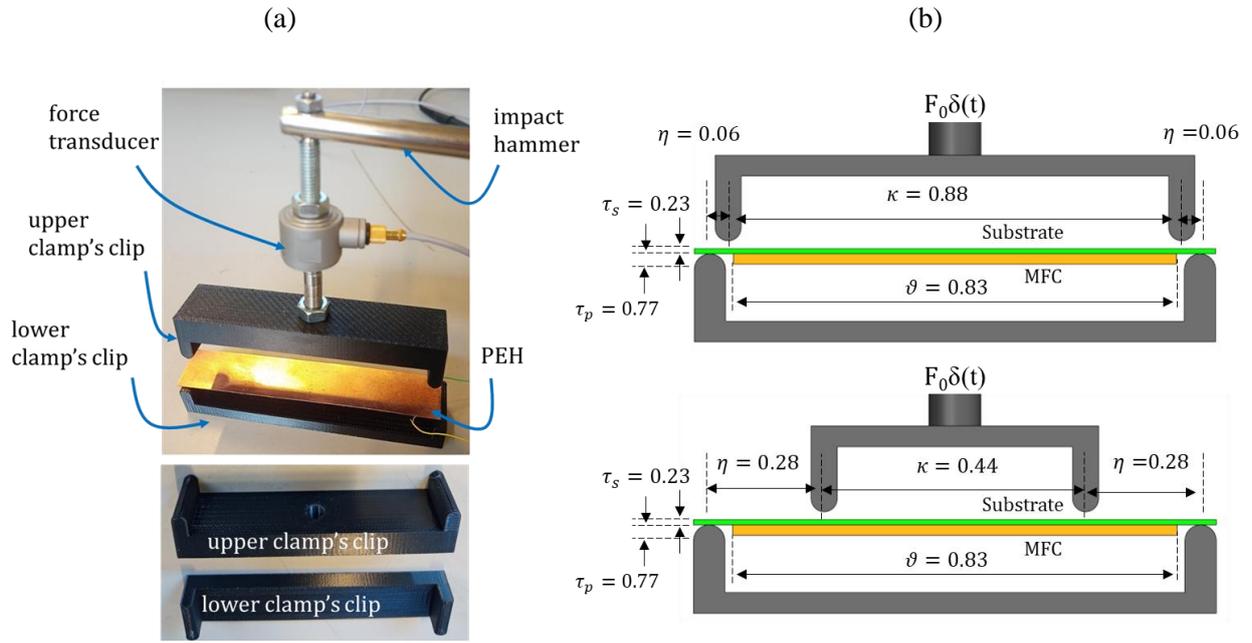

Figure 6. (a) Experimental hammer setup for model verification and (b) two FPB cases with different force applying spans with the dimensionless parameters

Specific experimental errors are common in this type of hammer impact. The exact impact locations are set visually and thus exposed to human error, as there might be slight differences in the intended and actual areas. Second, since there are two lines of contact on the upper clip for applying the load, if the upper clamp clip is not horizontal, the impacts from these two lines might not be transferred simultaneously. Third, the energy harvester weight is less than 7 g, so there is a possibility that the harvester slips on the clamp or has in-plane motion. Lastly, even thin wires, the electrical connection wires can interfere with the harvester's free vibration after impact. The first and second errors can be diminished by careful test operator and modal hammer skills. The third and fourth errors are unavoidable, and their contributions shall be carefully tracked; however, by connecting thin wires, the wire effects can be reduced.

The following material properties are employed in the modeling: The Young's moduli are 24.8 and 111.4 *GPa* for the piezoelectric and substrate, respectively. The corresponding densities are 5540 and 8960 *kg/m³*. The damping coefficient is derived from the author's previous experimental study on the same energy harvester sample [24] and considers frictional energy dissipation [25]. The relative dielectric coefficient and piezoelectric constant $d_{31}$ are 1800 and $-170\times10^{-12}$ *C/N,* respectively.

For this energy harvester, the undamped natural frequencies from the current model are compared with COMSOL software and a high-order Finite Element model [26] in Table 1. The natural frequencies for the four bending modes agree with the error below 5%. Significantly, the first natural frequency has less than 1% error compared with the COMSOL software result.

Table 1. The comparison of undamped natural frequencies between the presented method, COMSOL software, and high-order Finite Element model [26]

|  | Undamped natural frequencies (Hz) | | | |
| --- | --- | --- | --- | --- |
|  | COMSOL software | Finite element numerical model | The current method (presented in subsection 2.2) | Error (%) |
| First bending mode | 48.7 | 49.6 | 48.4 | 0.6 |
| Second bending mode | 196.7 | 198.7 | 193.4 | 1.7 |
| Third bending mode | 448.4 | 448.4 | 435.2 | 2.9 |
| Fourth bending mode | 807.7 | 799.7 | 773.7 | 4.2 |

Hammer impacts were varied from small to significant impact magnitudes to study the piezoelectric output voltage resulting from the input force. Figure 7 (a) shows a typical force measurement for κ=0.88 (*L*=90 mm) span. The impact duration is approximately 0.05 s, and the voltage generation peaks after the impact. The force-time plot shows that the hammer impact force is a single hit, which corresponds to the analytical impact assumption. Figure 7 (b) shows the impact force level at different tests κ=0.88 and κ=0.44. For κ=0.88, if the upper clamp span is close to the supports, the supports' responses are more severe; therefore, the measured force on the hammer becomes larger. That is why the measured force for the κ=0.88 span is 1-8 N while for κ=0.44, the range is 0.2-1.1 N.

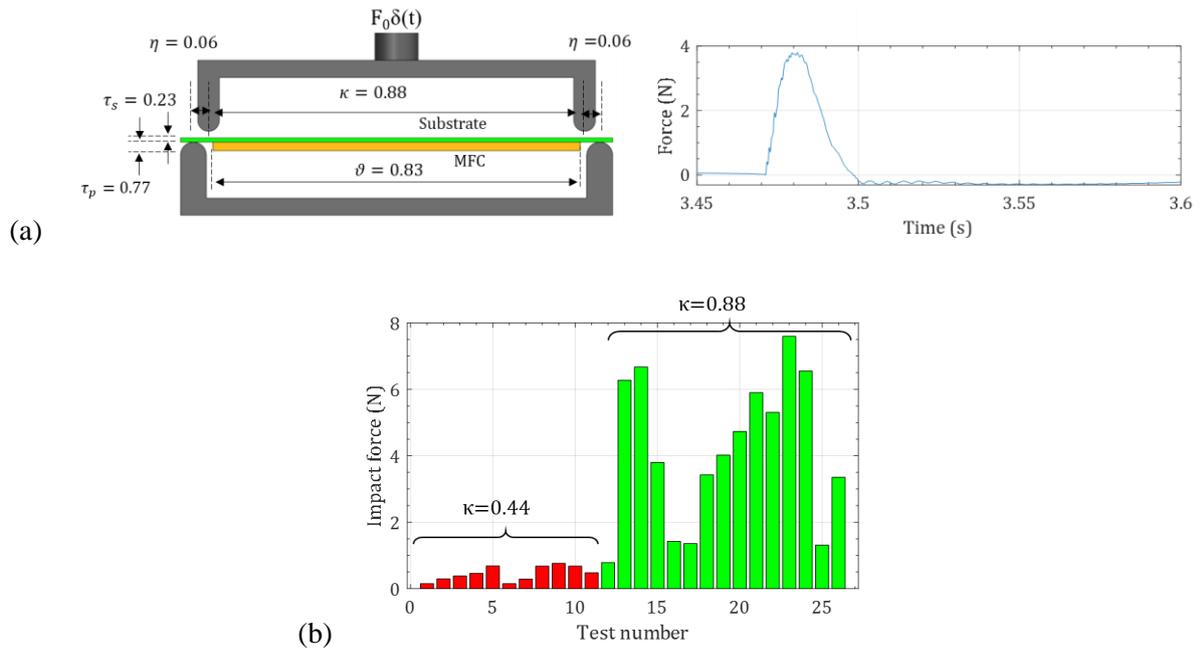

Figure 7. (a) A typical force and piezoelectric voltage measurements at κ=0.88 (*L*=90 mm) span, and (b) various impact force levels for κ=0.88 (*L*=90 mm) and κ=0.44 (*L*=45 mm).

Figure 8 compares the model and experimental voltage data for κ=0.44 and κ=0.88 spans. Two trends are discernible: the linear voltage generation increases with force impact amplitude, and the span (κ) has the effect that smaller κ leads to higher voltages. Both trends agree with the presented model. The experimental data confirm that the peak generated voltage increases approximately linearly with impact magnitude. The agreement between the model and data is reasonable for both κ=0.44 and κ=0.88 spans. In addition, the voltage generation for κ=0.44 is higher than for the κ=0.88 case because of the higher bending moment. This feature, further elaborated in subsection 3.4, validates the model.

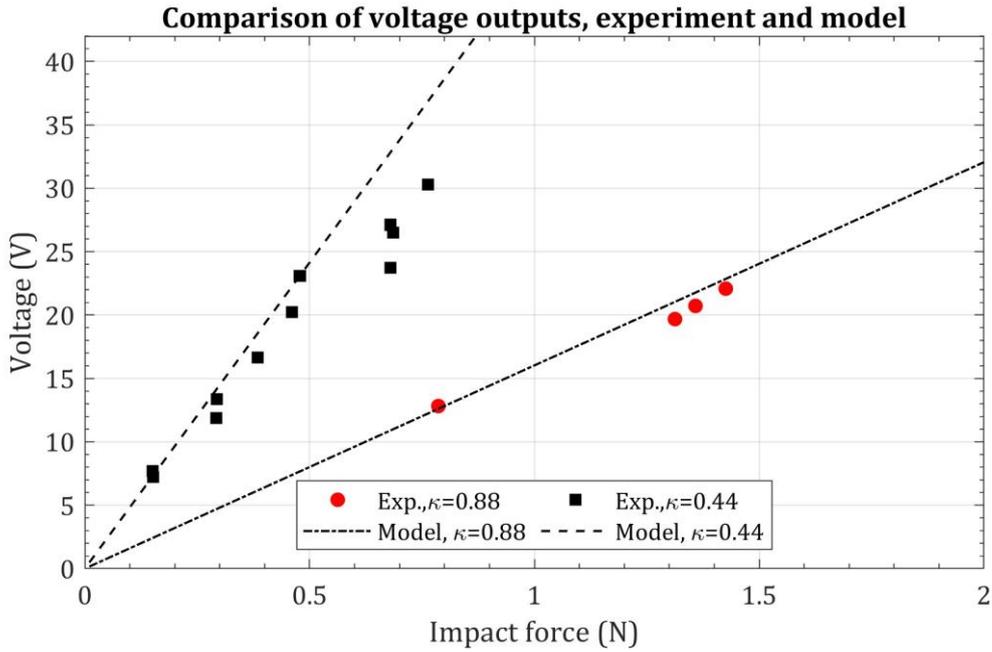

Figure 8. The comparison of FPB modeling results and experimental data under the impact with κ=0.44 and κ=0.88 spans.

### 3.3. Strain distribution in the FPB energy harvester

The strain distribution is a metric for the voltage generated in the PEHs and can be employed to optimize the energy generation over the PEH volume. The strain distribution is dictated by the boundary conditions imposed on the PEH. Figure 9 shows strain contours through the thickness and over the length of the harvester. As expected, the regions above and below the neutral axis have different signs since one area is under compression (negative strain), and one is under tension (positive strain). In addition, the greater distance to the neutral axis leads to greater strain magnitude; therefore, the axial strain at the MFC surface (piezoelectric layer) is maximum because of the MFC's surface location and thickness.

The mode-shapes impose onion-like layers in the strain contours (influenced by the boundary condition). These layered patterns are defined by the mode shapes, as depicted in Figure 9, especially by the first mode shape, which is the most significant mode. Figure 9 shows the strain patterns for three bending modes. Moreover, high strain values are found over most of the length of the PEH, ξ=10% to ξ=90%, and over 50% of the harvester thickness, thus occupying a considerable volume fraction of the PEH. Quantitatively, 48% of the piezoelectric volume is under at least 40% $\max(\varepsilon_{xx})$; showing that in the FPB a significant part of the piezoelectric layer contributes to the energy generation. Note that the piezoelectric volume contribution to the energy generation is approximately 33% for the cantilevered boundary condition.

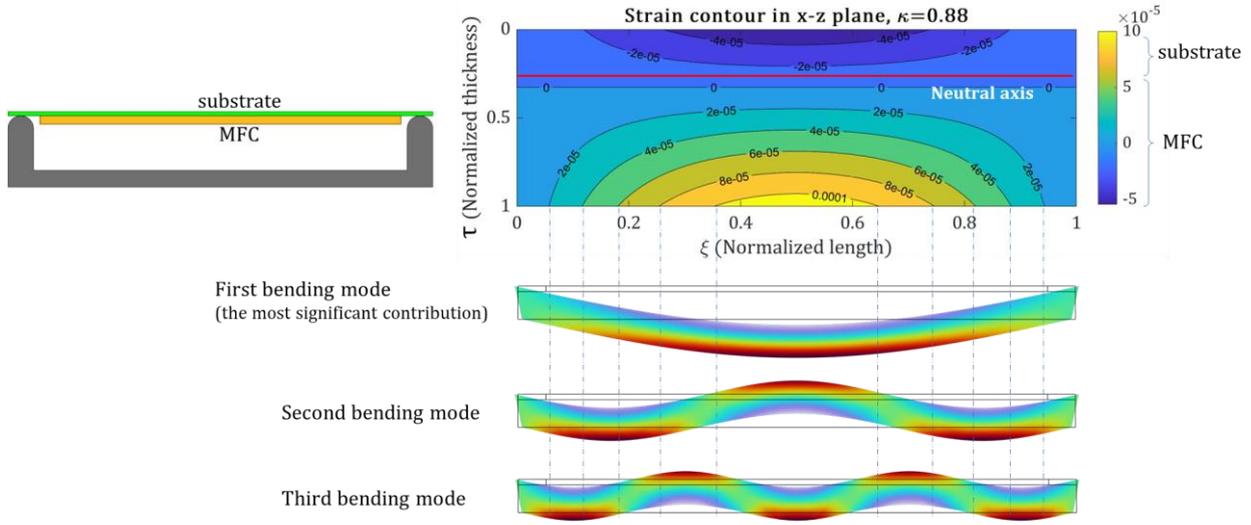

Figure 9. Strain contour in the reference FPB energy harvester over the length-thickness plane

### 3.4. Power generation and optimum electrical resistance load

Voltage and current peak outputs for the reference FPB harvester are investigated by the presented model with different resistance load connections and plotted in Figure 10 (a). As expected, voltage and current have converse trends. The power generation versus resistance load is plotted in Figure 10 (b). An optimum load leads to the highest power generation, similar to the typical power-load plots obtained for other piezoelectric energy harvesting systems. The maximum peak power under a unit impact was 6.7 mW/$N^2$ with the optimum load of 29kΩ.

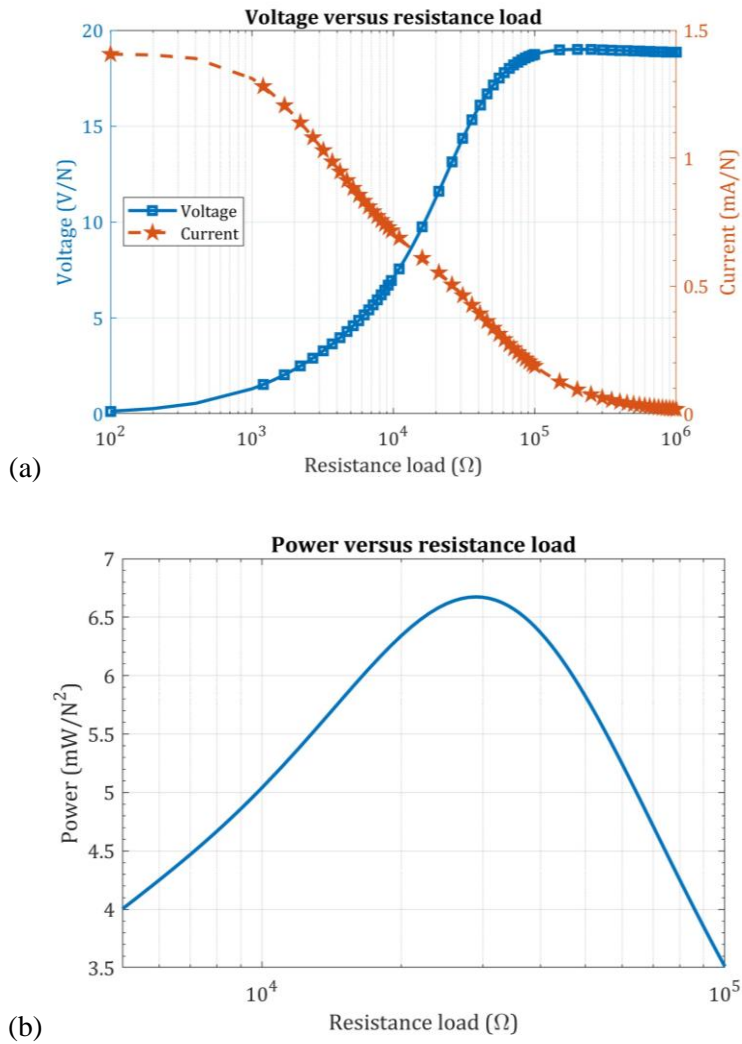

(a)

(b)

Figure 10. (a) voltage and current outputs and (b) power output versus load for the reference FPB harvester.

### 3.5. Sensitivity analysis of the FPB energy harvester

This section provides the sensitivity analysis of the FPB boundary condition parameters. Results are from the dynamic numerical solution of the analytical model for an impact force. Two dimensionless parameters are investigated, dimensionless force span $\kappa$ and piezoelectric length $L_p$ parameters.

### 3.5.1. Effect of force span (*L* or dimensionless parameter $\kappa$)

One central aspect of this FPB energy harvester is the top clamp span which affects the load transfer into the energy harvester and imposes notable change to the harvester dynamic. Thus, the effect of force span is studied here by varying the dimensionless parameter κ.

Figure 11 (a) and (b) compare the shear force and bending moment in the reference FPB harvester for the κ=0.88 (reference harvester) and κ=0.44 cases under a static load. The bending moment in FPB from point B to point C is constant, depending on the force magnitude and acting distance, e.g., A-B distance (η). Since $\eta_{\kappa=0.44} > \eta_{\kappa=0.88}$, the bending moment from B to C is greater for κ = 0.44. This is equally valid for the dynamic impact case. Recalling the modal force coefficient $\sigma_n$ from Eq. (17), $\sigma_n$ links to the mode shape at force contact points, e.g., $\phi_n|_B$ and $\phi_n|_C$. Therefore, force contact points close to end-clamps have smaller modal force coefficient, leading to smaller force transfer to the structure.

$$\sigma_n = -\left(\left(1 - \frac{U_1}{L}\right) \cdot \phi_n|_B + \frac{U_1}{L} \cdot \phi_n|_C\right) \qquad \text{Eq. (17)}$$

Thus, reducing κ will increase the bending moment and the strain on the PEH. This is shown in Figure 11 (c), where the strain in the x-z plane is plotted for the κ = 0.44 and κ = 0.88 cases under a unit impact force. The strain contours have the same patterns; however, the strain amplitude for the κ = 0.44 case is considerably greater than the κ = 0.88 case. Like the strain, the beam deformation is influenced by κ. The maximum deflection with κ = 0.44 is 1.61mm, which is approximately four times greater than 0.41mm for κ = 0.88. The more significant deformation and strain for smaller $\kappa$ is due to the modal force coefficient given in Eq. (17). The variation of deformation and strain with force span is nonlinear, as the mode-shapes are nonlinear.

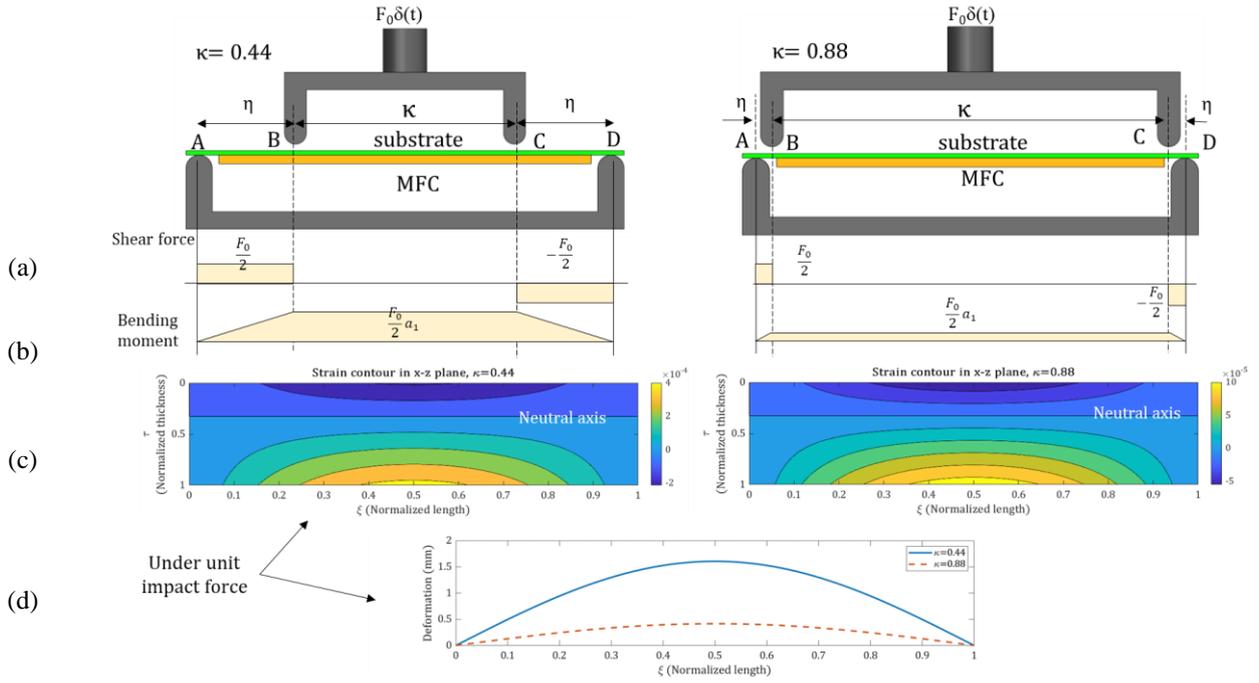

Figure 11. Comparing (a) shear force, (b) bending moment, (c) axial strain in the FPB, and (d) maximum beam deformation of the reference FPB harvester with κ = 0.44 and κ = 0.88 spans.

If one is interested in the voltage generation per unit strain, the $\lambda$ index can be compared between these two spans, as given in Table 2. Greater $\lambda$ implies better use of material for energy generation. As presented in Table 2, $\lambda_{\kappa=0.88}$ is 25% greater than $\lambda_{\kappa=0.44}$, demonstrating that a better energy conversion performance is obtained for greater $\kappa$ though the peak voltage is lower. These opposite trends indicate that the best energy conversion mechanism differs according to the objective of voltage generation or energy conversion, and both aspects should be considered by designers.

Table 2. The PEH performance for κ = 0.44 and κ = 0.88 cases.

| Performance under the unit impact force | $\max(w)$ | $\max(\varepsilon_{xx})$ | Peak voltage $V_p$ (V) | $\lambda = \dfrac{V_p}{\max(\varepsilon_{xx})}$ |
|---|---|---|---|---|
| κ = 0.44 | $1.61\ mm$ | $4.32\times10^{-4}$ | $57.79\ V$ | $1.34\times10^5\ V$ |
| κ = 0.88 | $0.41\ mm$ | $1.12\times10^{-4}$ | $18.81\ V$ | $1.68\times10^5\ V$ |

By investigating the variation of the $\lambda = \dfrac{V_p}{\max(\varepsilon_{xx})}$ coefficient with the non-dimensional load span (κ), the best energy conversion design can be demonstrated. $\lambda$ is plotted against the span κ in Figure 12. By increasing κ up to 0.83, the $\lambda$ coefficient increases, indicating that the energy harvester will generate a higher voltage per unit strain. However, for κ > 0.83, $\lambda$ shows a slight decline.

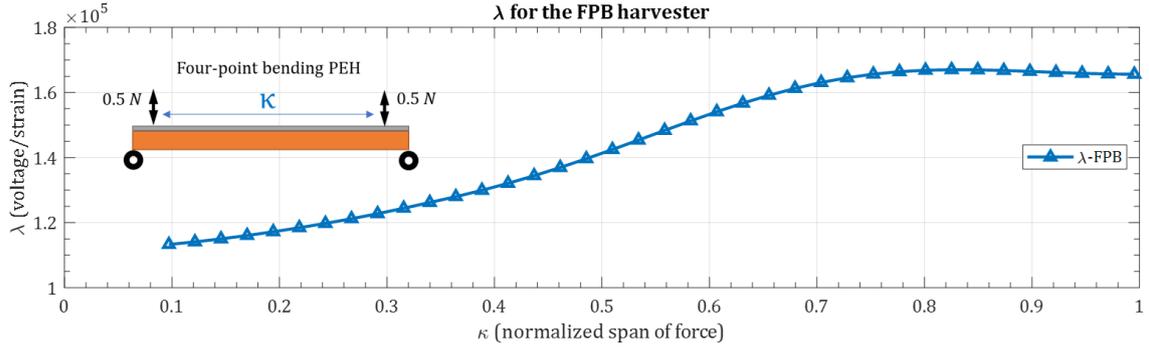

Figure 12. The $\lambda = \frac{\text{Voltage}}{\text{Axial strain}}$ coefficient versus the span of applying force for the reference FPB harvester.

### 3.5.2. Effect of Piezoelectric layer length ( $L_p$ or dimensionless parameter $\vartheta$)

Figure 13 (a) shows maximum power density under a unit impact for five different piezoelectric layer length values, characterized by $\vartheta$= 0.75, 0.79, 0.83, 0.86, and 0.9. Note that the MFC-piezoelectric layer is symmetrically centered, and $\vartheta$ indicates the fraction of the harvester length occupied by the piezoelectric layer. Increasing the piezoelectric layer length means more piezoelectric material, generating more power. This expected power density increase is obvious in Figure 13 (a); increasing piezoelectric length increases the power density. In addition, the optimum load is slightly reduced by the piezoelectric length increase.

A piezoelectric length increase will enhance the power density, as expected. However, this increase is not linear because of the nonlinear strain in the FPB harvester, as the 1st bending mode is shown in Figure 13 (a). The power density variation due to the piezoelectric size is also given in Figure 13 (b); overall, increasing piezoelectric length by 3 to 4% increases the power density by over 10%. In addition, increasing the piezoelectric size changes the electromechanical conversion efficiency, leading to a slight reduction of the optimum load.

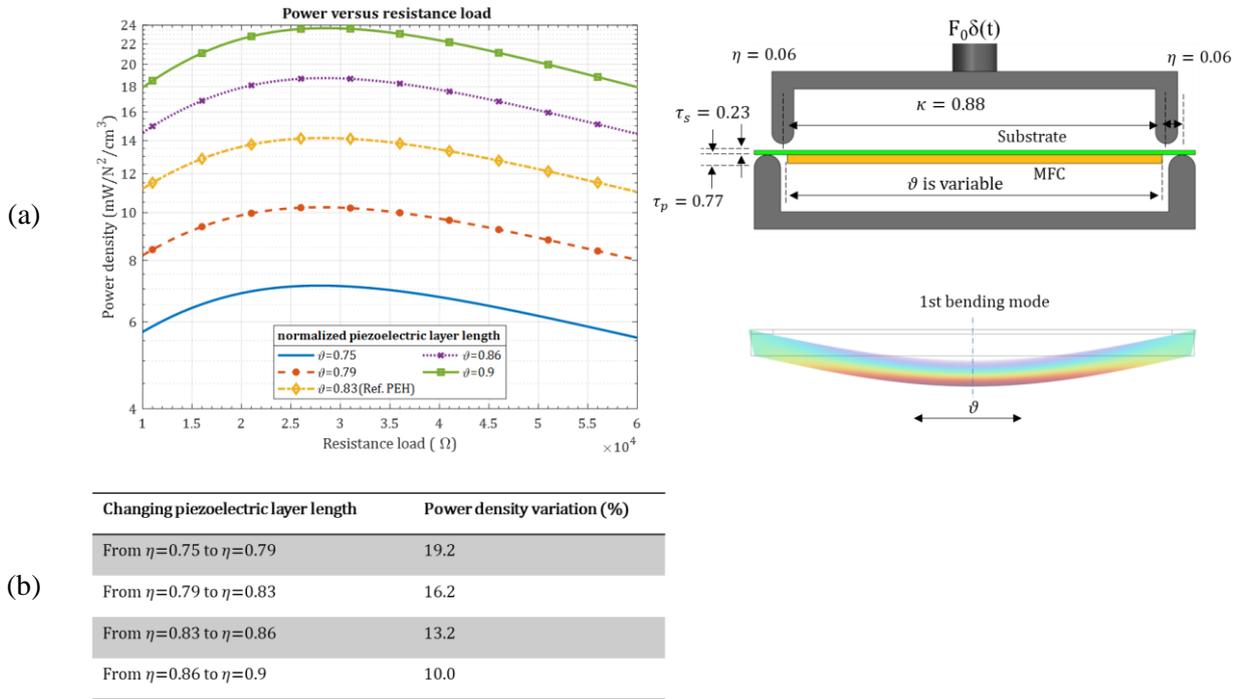

Figure 13. (a) The optimum resistance load and power density for different normalized piezoelectric length values, $\vartheta$ = 0.75, 0.79, 0.83, 0.86, and 0.9, (b) power density versus $\vartheta$, and (c) power density variation for different normalized piezoelectric length values.

### 3.6. Comparison of four-point bending and cantilevered harvesters

This section compares the FPB energy harvester with the typical cantilevered harvesters concerning the mechanical stress and the harvester voltage generation. Results are calculated from the presented model under an impact force.

Figure 14 compares the performance of FPB and cantilevered harvesters; Figure 14 (a) indicates the voltage and maximum stress on the piezoelectric material, and Figure 14 (b) presents the $\lambda$ coefficient. The x-axis variable is the force span dimensionless parameter κ.

From Figure 14 (a), in the FPB boundary condition, increasing the normalized force span, κ, reduces the force transmission to the energy harvester (applied load moves closer to supports), as described in subsection 3.4. This phenomenon reduces the mechanical stress on the piezoelectric layers and leads to voltage reduction. Note that for the extreme case of κ = 1, the inner and outer clamps reach each other; thus, force is not transferred to the harvester, resulting in zero voltage generation. On the other hand, in the cantilevered beam, increasing the force span increases the stress on the piezoelectric material (load moves away from support), leading to higher voltage generation.

From Figure 14 (a), the maximum stress at the highest voltage generation in the cantilevered is 32.7 MPa/N, which is considerably higher than 12.9 MPa/N for the FPB harvester. Typical cantilevered harvesters under the base excitations work at the highest stress point. Moreover, the stress in the cantilevered harvester dramatically soars by changing the force span. Regarding the voltage generation, the maximum cantilevered voltage output is 13.5% higher than the FPB voltage output, though this increase comes with a price of withstanding a 250% higher mechanical stress. Overall, the power generation performance for the FPB harvester is comparable with that of the cantilevered harvester, while the FPB boundary conditions result in substantially reduced stress.

As an energy conversion per unit strain metric, the $\lambda$ coefficient is plotted in Figure 14 (b) against the force span. Overall, $\lambda$ for the FPB harvester is substantially higher than for the cantilevered harvester, demonstrating a better energy conversion in the FPB device. Moreover, $\lambda$ is linked to the force span, increasing with span in the FPB harvester, while $\lambda$ reduces with span increase for the cantilevered harvester. The cantilevered harvester undergoes smaller strain by reducing span. From the $\lambda$ metric, the FPB harvester with $\kappa=0.83$ force span appears optimal. The best energy conversion coefficient is $\lambda=1.67\times10^5$ V/$\varepsilon$ for the FPB harvester, which is 315% higher than the typical cantilevered harvester $\lambda=0.54\times10^5$ V/$\varepsilon$.

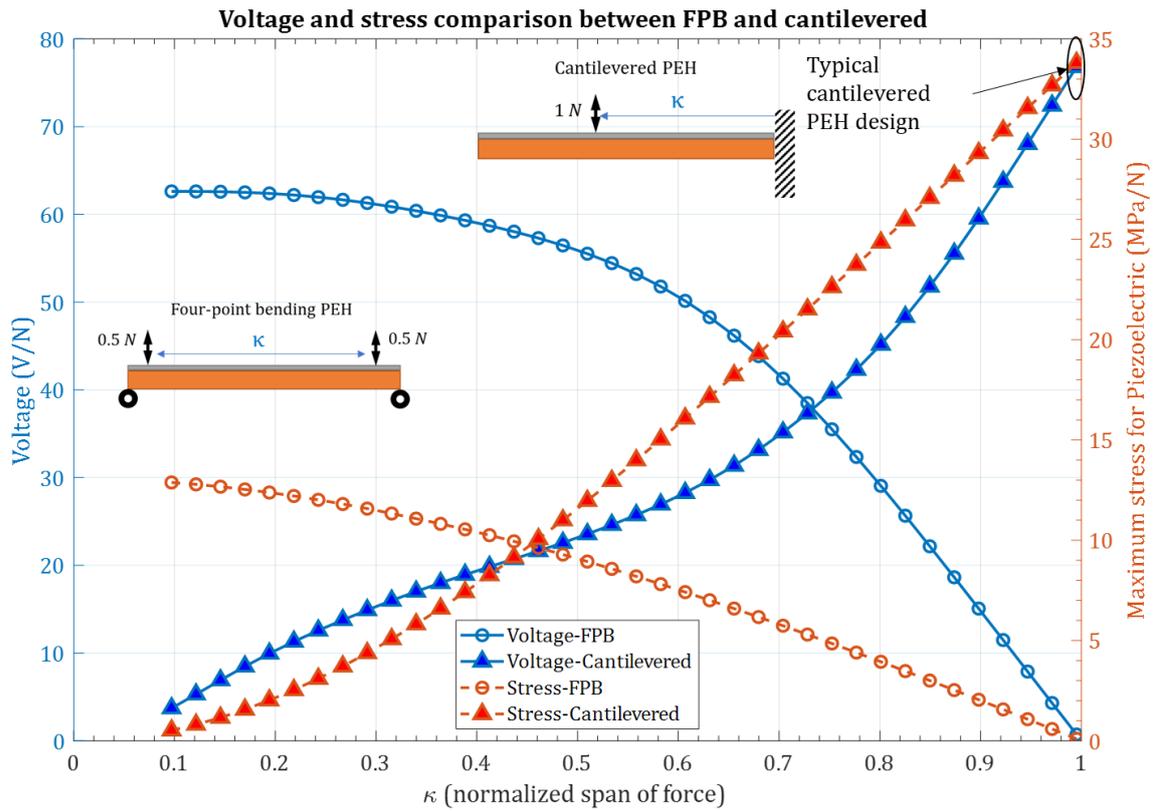

(a)

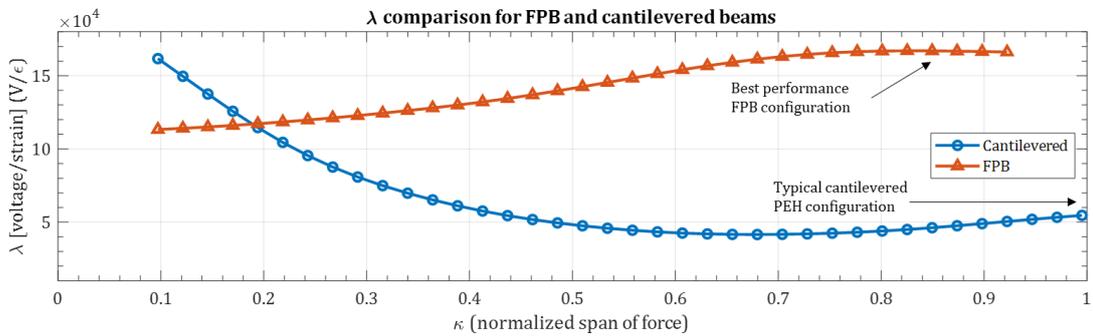

(b)

Figure 14. The comparison between the four-point bending (FPB) and cantilevered harvesters versus the span of applying load, (a) voltage and maximum stress for the Piezoelectric layer, and (b) the $\lambda = \frac{\text{Voltage}}{\text{Axial strain}}$ coefficient.

### 3.7. FPB conversion performance and Figure of Merit

The overall performance of the FPB-PEH is demonstrated in Figure 15 by plotting Voltage, $\lambda$ (the energy conversion index), and $FoM$ (Figure of Merit) versus the normalized force span. The best voltage generation with κ=0.1 corresponds to the lowest energy conversion performance, meaning that a substantial part of the mechanical energy is wasted. On the other hand, the best energy conversion (the highest $\lambda$) point generates small voltage compared to the maximum available voltage.

Considering these two parameters, the best Figure of Merit ($FoM$) corresponds to moderate energy conversion performance and voltage generation. The maximum $FoM$ occurs at κ=0.53, giving the energy conversion index of $\lambda=1.45\times10^5\ V/\epsilon$ and $V_{OC}$=54.4 V/N.

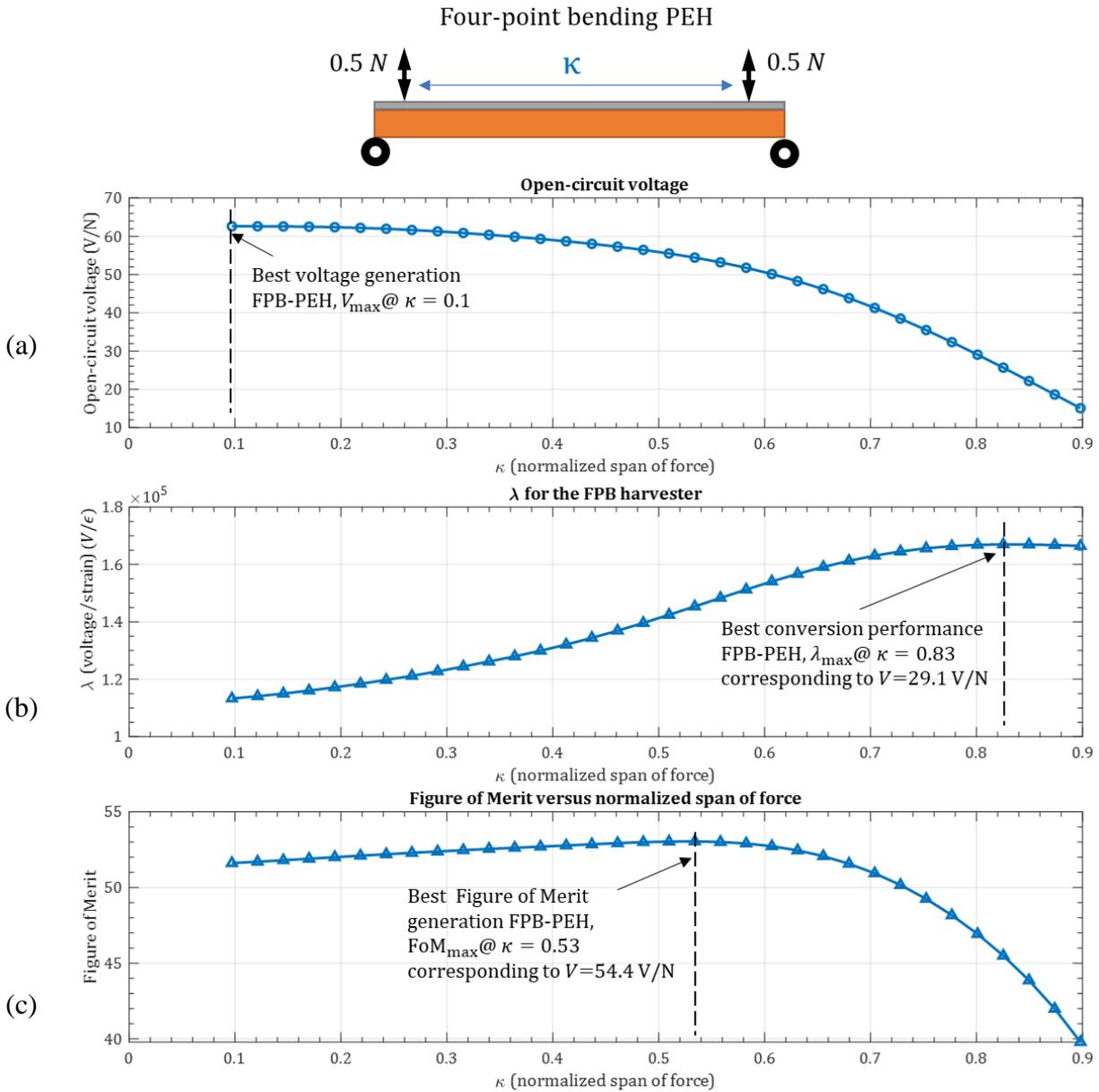

Figure 15. The overall performance of the presented four-point bending (FPB)-PEH in terms of (a) voltage generation, (b) energy conversion performance ($\lambda = \frac{\text{Voltage}}{\text{Axial strain}}$), and (c) Figure of Merit.

The energy harvester's reliability is significant since autonomous self-powering systems require reliable long-term performance. Therefore, the harvester's fatigue life plays a vital role in the harvester's design. The reliable upper limit of the dynamic strain for the MFC is estimated to be 600 με, and above 1000 με cracks will initiate [27]. The accelerated life fatigue tests have demonstrated that reducing the applied stress by 21% can double an MFC's cyclic life [28]. This study proposes a novel four-point-bend boundary condition that has reduced the applied strain. Thus, the presented FPB harvester may increase the cyclic life of piezoelectric energy harvesters. Figure 16 compares the maximum strain for the FPB harvester and a typical cantilevered harvester. As shown in Figure 16, the strain in the FPB harvester with the best Figure of Merit is $\epsilon_{xx}$=374.5 με/N, while a typical PEH cantilever design has a strain value of above $\epsilon_{xx}$=1214.9 με/N. Therefore, while the FPB harvester can withstand a safe cyclic life with 1.6 N load, a typical cantilevered harvester will be damaged by 1 N load. The strain in the piezoelectric material can be significantly reduced by this innovative boundary condition, while maintaining performance. Hence this will significantly improve the reliability of piezoelectric energy harvesters.

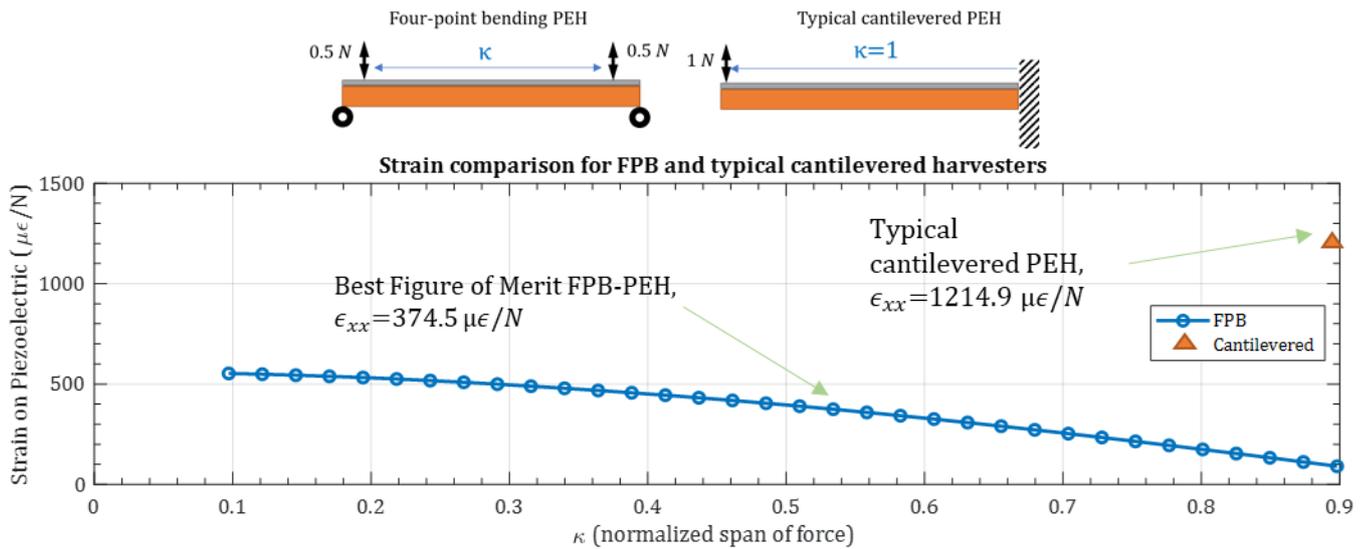

Figure 16. A comparison of the strain between the proposed FPB and cantilevered energy harvesters.

## 4. Concluding remarks and future works

This paper proposes a new piezoelectric energy harvester (PEH) design with a better energy conversion performance. The four-point-bend piezoelectric energy harvester (FPB-PEH) distributes strain more uniformly than typical cantilevered harvesters. This results in a reduced strain magnitude, demonstrating an improved use of the material for energy conversion. An electromechanical coupled model of the FPB-PEH and analytical expressions for the energy harvester's voltage, mechanical deformation, and strain were obtained. The model has been verified against experimental data and Finite Element results. A reference FPB-PEH with one Macro Fiber Composite (MFC) piezoelectric layer was studied.

Strain fields over the PEH volume were analyzed, and it is concluded that 48% of the piezoelectric contributes significantly to energy harvesting, which is a considerable improvement relative to a cantilevered PEH. An energy conversion coefficient per unit strain, λ, was used as a performance metric. Comparisons between the performance of the FPB and cantilevered PEHs were made, and the superiority of the proposed FPB was identified in energy conversion per unit strain and energy generation uniformity. The performance metric λ was studied for the FPB and cantilevered PEHs under different applied load spans. Sensitivity analyses for the piezoelectric layer length, the piezoelectric layer's location, and the applied force span were presented. This paper contributes to the development of better and more reliable piezoelectric energy harvesters and innovative design for impact-based energy harvesters. More experimental investigationof the contact regions and clamp types are proposed for future work. Moreover, further studies on the parameters of the energy harvesters, such as thickness and multi-layered designs, and various material investigations are needed.

## Declaration of competing interests

The authors declare that they have no known competing interests that may influence the scientific works in this paper.

## Acknowledgments

This research is partially financed by the Independent Research Fund Denmark International Post-doc grant under grant number 1031-00001B and the Lundbeck LF-Experiment grant.